\documentclass[12pt]{amsart} 
\usepackage{amsmath,amsfonts,amssymb,amsthm}

\pagespan{1}{?}

\usepackage{epsfig}
\usepackage{graphicx}

\theoremstyle{definition}

\numberwithin{equation}{section}

\newcommand\Lie[1]{{\mathfrak{#1}}}
\def\R{{\mathbb R}}

\begin{document}
\title[A Liouville type theorem for Carnot groups]
{A Liouville type theorem for Carnot groups: a case study }

\author[Alessandro Ottazzi]{\noindent Alessandro Ottazzi}
\email{alessandro.ottazzi@unimib.it}
\address{\newline Department of Mathematics,
\newline University of Milano--Bicocca
\newline Milano , 20126 ,
\newline Italy}



\author[Ben Warhurst]{\noindent Ben Warhurst}
 \email{benwarhurst68@gmail.com}
 \address{\newline Department of Mathematics,
 \newline Milano , 20126 ,
 \newline University of Milano--Bicocca, Italy}



\begin{abstract}
In~\cite{CC}, the authors show that if $\phi$ is  1-quasiconformal  on an open subset of a Carnot group G,  then composition with $\phi$ preserves Q-harmonic functions, where Q denotes the homogeneous dimension of G. Then they combine this with a regularity theorem for Q-harmonic functions to show that $\phi$ is in fact $C^\infty$. As an application, they observe that a Liouville type theorem holds for some Carnot groups of step 2.

In this article we argue, using the Engel group as an example, that a Liouville type theorem can be proved for every Carnot group. Indeed, the fact that 1-quasiconformal maps are smooth allows us to obtain a Liouville type theorem by applying the Tanaka prolongation theory~\cite{tanak1}.

\end{abstract}

\maketitle


\section{Introduction}

The classical Liouville theorem states that $C^4$--conformal maps between domains of $\R^3$ are the restriction of the action of some element of the group $O(1,4)$. The same result holds in $\R^n$ when $n> 3$ (see, e.g., Nevanlinna~\cite{N}). A major advance in the theory was the passage from smoothness assumptions to metric assumptions (see Gehring~\cite{Geh} and Reshetnyak~\cite{Res}): the conclusion of Liouville's theorem holds for 1-quasiconformal maps. When the ambient space is not Riemannian there are similar theorems. Capogna and Cowling proved in~\cite{CC} that 1-quasiconformal maps defined on open subsets of a Carnot group $G$ are smooth and applied this result to give some Liouville type results. In particular, it is now known that 1-quasiconformal maps between open subsets of H-type groups whose Lie algebra has dimension larger than $2$ form a finite dimensional space. This follows by combining  the smoothness result in~\cite{CC} and the work of Reimann~\cite{R}, who established the corresponding infinitesimal result. Moreover, if $G$ is a Carnot group of step two such that the strata preserving automorphisms are all dilations, then 1-quasiconformal maps are translations composed with dilations~\cite{CC}.

In this article we combine the smoothness result in~\cite{CC} with the Tanaka prolongation theory to show that when $G$ is the Engel group (step three), 1-quasiconformal maps form a finite-dimensional space. A characterization of conformal maps for the Engel group can also be found in~\cite{cdkr0}. There the Engel group is viewed as the nilradical of the group ${\rm Sp}(2,\R)$ and the theory of semisimple Lie groups plays a central role.
The main point of interest in our approach is that it extends to all Carnot groups.  We shall treat the general case in a forthcoming paper.

The next section is devoted to the study of 1-quasiconformal maps on the Engel group. First we define the basic formalism: we introduce the contact structure, the subriemannian metric and we give the definition of quasiconformal maps. In particular, we interpret 1-quasiconformal maps as conformal transformations. Next, we restrict the attention  to  the infinitesimal level introducing conformal vector fields. This leads to a system of differential equations. We then proceed by constructing a prolongation of these differential equations formalized in terms of Tanaka prolongation. The latter is a graded Lie algebra that ends up to be isomorphic to the Lie algebra of conformal vector fields. Finally, we use a standard argument to show that any conformal map is the restriction of the action of some element in the automorphism group of the prolongation algebra.

 \section{A case study: the Engel group}
\subsection{Notation and definitions.}
Let $\Lie{g}$ be the real Lie algebra generated by the vectors $X_1, X_2, Y, Z$ and the nonzero brackets  $[X_1,X_2]=Y$ and $[X_1,Y]=Z$.
This is a stratified nilpotent Lie algebra of step three. Namely $\Lie{g}=\Lie{g}_{-1}+\Lie{g}_{-2}+\Lie{g}_{-3}$, where $\Lie{g}_{-1}={\rm span}\{X_1,X_2\}$, $\Lie{g}_{-2}=\R Y$ and $\Lie{g}_{-3}=\R Z$. The symbol $+$ will denote the direct sum of vector spaces.
We write $G$ for the connected and simply connected Lie group whose Lie algebra is
$\Lie{g}$. We choose exponential coordinates  $(x_1,x_2,y,z)={\rm exp}(x_2X_2 +yY +zZ){\rm exp}(x_1X_1)$. We identify the Lie algebra with the tangent space $T_e G$ to $G$ at the identity $e$, and for $X$ in $\Lie{g}$ we write $\tilde{X}$ for the left--invariant vector field that agrees with $X$ at $e$. The left invariant vector fields corresponding to the basis vectors are
\begin{align*}
\tilde{X}_1 &= \frac{\partial}{\partial x_1}\qquad \tilde{X}_2 = \frac{\partial}{\partial x_2}+x_1\frac{\partial}{\partial y}+\frac{x_1^2}{2}\frac{\partial}{\partial z}\\
\tilde{Y}&=\frac{\partial}{\partial y}+x_1\frac{\partial}{\partial z}\\
\tilde{Z}&=\frac{\partial}{\partial z}
.\end{align*}
 The vector fields in $\tilde{\Lie{g}}_{-1}$ determine at each point a subspace of the tangent space that we call {\it horizontal space}. These subspaces vary smoothly from point to point and give rise to the horizontal tangent bundle. Since $\Lie{g}_{-1}$ generates $\Lie{g}$ with the brackets, the sections of the horizontal tangent bundle generate all possible vector fields with linear combinations of commutators. We define an inner product $\langle \,,\, \rangle$ on $\Lie{g}$ for which the different layers are orthogonal and we transport the inner product to the tangent space at each point using left translation. We denote by $\langle \,,\, \rangle_p$ the inner product in $T_pG$, where $p\in G$.
 This allows us to define a left-invariant Carnot--Carath\'eodory metric $d$ on $G$, as follows. A smooth curve is said to be {\it horizontal} if its tangent vectors are horizontal. The length of a horizontal curve is the integral of the lengths of its tangent vectors. The distance between two points is then the infimum of the lengths of the horizontal curves joining them. A nilpotent stratified Lie group  with such a metric is called Carnot group and the example $G$ we are considering is known as the Engel group.

A diffeomorphism $\phi$ between open sets of $G$ is called a {\it contact mapping} if its differential $\phi_*$ preserves at each point the horizontal space. Given a 1-parameter group of contact mappings, we denote by $V$ the corresponding infinitesimal generator and call it a {\it contact vector field}.
Contact vector fields are characterized by the differential equations arising from the condition that for every horizontal vector field $\tilde{X}$ one has $[V,\tilde{X}]= f\tilde{X}_1
+ g\tilde{X}_2$ for some functions $f$ and $g$.
It is well known~\cite{War1, otz1} that for  the Engel group the space of contact vector fields is infinite dimensional. Indeed, it is straightforward to prove that $V= f\tilde{Z}+ \tilde{X}_1 f \tilde{Y}+
\tilde{X}_1^2 f\tilde{X}_2$ is a contact vector field for every smooth function $f=f(x_1)$.

Let $\mathcal{U},\mathcal{V}\subset G$ be open domains, and let $\phi :\mathcal{U}\rightarrow \mathcal{V}$  be a homeomorphism. For $p\in \mathcal{U}$ and for small $t\in \R$ we define the distortion as
$$
H_{\phi}(p,t)= \frac{{\rm max}\{d(\phi(p),\phi(q)) | d(p,q)=t\}}{{\rm min}\{d(\phi(p),\phi(q))|d(p,q)=t\}}.
$$
We say that $\phi$ is {\it quasiconformal} if there exists a constant $\lambda$ such that
$$
\lim \sup_{t\rightarrow 0} H_\phi (p,t)\leq \lambda
$$
 for all $p\in \mathcal{U}$. Furthermore, $\phi$ is locally quasiconformal if it is quasiconformal in a neighborhood of each point and 1-quasiconformal when $\lambda =1$.

 In~\cite{CC}, the authors prove that 1-quasiconformal mappings are smooth and characterized by the conditions of being locally quasiconformal with Pansu differential  $D\phi(p) $ coinciding with a similarity at every point $p$   (that is, a product of a dilation and an isometry). In fact, $\phi$ is contact~\cite{Pu90}. We remind the reader that
$$
D\phi (p) \,q= \lim_{t\rightarrow 0} \delta_t^{-1}[\phi(p)^{-1}\phi(p\delta_tq)],
$$
where $p,q\in G$ and for every $t\in \R^+$, $\delta_t$ denotes the automorphic dilation that in this case study is $\delta_t(x_1,x_2,y,z)=(e^tx_1,e^t x_2,e^{2t}y,e^{3t}z)$. The Pansu differential is an automorphism of $G$ and by means of the Baker-Campbell-Hausdorff formula defines an automorphism of the Lie algebra $\Lie{g}$ as well, that we denote by $d\phi (p)$. The condition that $D\phi(p) $ coincides with a similarity  implies that $d\phi (p)$  is a similarity when restricted to $\Lie{g}_{-1}$~\cite[Lemma 5.2]{CC}. 
 Writing the differential $\phi_*$ at  $p$  according to the basis of left--invariant vector fields,  it turns out that $\phi_*$ and the Pansu differential are equal when restricted to the horizontal space. Therefore, the condition for 1-quasiconformality is equivalent to $\phi_*$ being a similarity of the horizontal space at every point. More explicitly, for all horizontal vectors $X,X^\prime \in T_pG$ we have
\begin{align*}
\langle \phi_{*p}^{tr}\phi_{*p}X,X^\prime \rangle_p&=\langle \phi_{*p}X,\phi_{*p}X^\prime \rangle_{\phi(p)}  \\
&=k\langle X,X^\prime \rangle_p
\end{align*}
whence for every $p\in\mathcal{U}$, the map $\phi_{*p}$ restricted to the horizontal space lies in the two dimensional conformal group $CO(2)= \{A \in GL(2,\R) | A A^{tr}=k I\}$.
In view of this fact we shall also refer to 1-quasiconformal maps as {\it conformal} maps. Furthermore, since the inner product is left invariant, we observe that left translations $l_p$ are conformal.

\subsection{Conformal vector fields.}
Let  $\phi_t$ be a 1--parameter group of conformal maps defined on an open set $\mathcal{U}$ and let $V$ be the vector field whose local flow is $\phi_t$. Then
\begin{align*}
\frac{d}{dt}(\phi_t)_* (\tilde{X})_p {|_{t=0}}&= \lim_{t\rightarrow 0} \frac{(\phi_{-t})_{*\phi_t(p)}\tilde{X}_{\phi_t(p)}-\tilde{X}_p}{t}\\
&= \lim_{t\rightarrow 0} \frac{(\phi_{-t})_{*\phi_t(p)}l_{{\phi_t(p)}_{*e}}-l_{{p}_{*e}}}{t} \tilde{X}_e\\
 &= [\tilde{X},V]_p =-{\rm ad}V (\tilde{X})_p.
\end{align*}
If $\tilde{X}$ is horizontal, we conclude from the chain of equalities above that $-{\rm ad}V$ restricted to the horizontal space is in  $\Lie{co}(2)=\{A\in\Lie{gl}(2,\R)|A+A^{tr}=kI\}$.

Using the basis of left invariant vector fields, we can write $V= f_1\tilde{X}_1 +f_2\tilde{X}_2+g\tilde{Y}+h\tilde{Z}$, where the coefficients are smooth functions. Then the contact conditions for $V$ are
\begin{align*}
[V,\tilde{X}_1]&=a\tilde{X}_1 +b\tilde{X}_2\\
[V,\tilde{X}_2]&=c\tilde{X}_1 +d\tilde{X}_2,
\end{align*}
for some functions $a,b,c,d$, which imply the following system of differential equations
\begin{equation}\label{contact}
\begin{cases}
\tilde{X}_1 g=-f_2\\ \tilde{X}_1h=-g
\end{cases}
\qquad\begin{cases}
\tilde{X}_2 g=f_1\\\tilde{X}_2h=0
\end{cases}
\qquad \tilde{Y}h=f_1.
\end{equation}
The condition of conformality on $V$ implies that the matrix
$$
{\rm ad}(V)_{|{\Lie{g}_{-1}}}(p)= \bmatrix \tilde{X}_1f_1 (p) & \tilde{X}_2 f_1(p)\\ \tilde{X}_1
f_2(p) & \tilde{X}_2 f_2 (p) \endbmatrix
$$
lies in $\Lie{co}(2)$, whence
\begin{equation}\label{conformal}
\begin{cases}
\tilde{X}_1f_1 &= \tilde{X}_2 f_2\\
\tilde{X}_2 f_1 &= -\tilde{X}_1 f_2.
\end{cases}
\end{equation}
A vector field that satisfies~\eqref{contact} and~\eqref{conformal} on its domain of definition is said to be a {\it conformal vector field}.

\subsection{Prolongation of the differential equations.}
In order to gather information on the space of functions that solve~\eqref{contact} and~\eqref{conformal} we consider higher order derivatives: roughly speaking, if the derivatives in all directions of the coefficients of $V$ vanish at a certain order, we may conclude that $V$ has polynomial coefficients and therefore varies in a finite dimensional space. For abelian differential operators, this idea  is formalized by the prolongation of Singer and Sternberg~\cite{ss}. In the noncommutative case, as is the situation we are interested in, the procedure was generalized by Tanaka~\cite{tanak1} and it was used to generalize the study of infinitesimal automorphisms of G-structures by different authors. For the contact structures, the Tanaka prolongation theory was used in~\cite{yam} and more recently by the authors of this article in different collaborations~\cite{alefil,OWnote,OW1,war2}.

Instead of giving the general argument of the Tanaka prolongation theory, we choose here to illustrate the method in our case study.
In order to do that, we first fix the following notations. For every $p\in \mathcal{U}$ we define
\begin{align*}
A^{-1}_V (p)&= (f_1(p),f_2(p),0,0)\\
A^{-2}_V (p)&= (0,0,g(p),0)\\
A^{-3}_V (p) &=(0,0,0,h(p)),
\end{align*}
where we interpret the vectors on the right hand side as elements of $\Lie{g}$ according to the fixed basis.
It follows that the conditions in ~\eqref{contact} are equivalent to
$$
[A^j_V(p),X]=\tilde{X}(A^{j-1}_V)(p),
$$
where $X\in \Lie{g}_j$ and $j=-1,-2$.

Next we define a graded homomorphism of  vector spaces $A^0_V : \Lie{g}\rightarrow \Lie{g}$ by setting
$$
A^0_V (p)(X) = \tilde{X}(A^j_V)(p),
$$
where $X\in \Lie{g}_j$, $j=-1,-2,-3$. Namely
\begin{align*}
A^0_V (p)(X_1) &= \tilde{X}_1 (A^{-1}_V)(p) = (\tilde{X}_1f_1(p),\tilde{X}_1f_2(p),0,0)\\
A^0_V(p) (X_2) &= \tilde{X}_2 (A^{-1}_V)(p) = (\tilde{X}_2f_1(p),\tilde{X}_2f_2(p),0,0)\\
A^0_V (p) (Y) &= \tilde{Y} (A^{-2}_V)(p)= (0,0,\tilde{Y}g(p),0)\\
A^0_V (p)(Z)&= \tilde{Z}(A^{-3}_V)(p)= (0,0,0,\tilde{Z}h(p)),
\end{align*}
or equivalently
\begin{equation}\label{A0}
A_V^0(p)=\bmatrix
\tilde{X}_1f_1(p)&\tilde{X}_2f_1(p)&0&0\\
\tilde{X}_1f_2(p)&\tilde{X}_2f_2(p)&0&0\\
0&0 &\tilde{Y}g(p)&0\\
0&0&0&\tilde{Z}h(p)
\endbmatrix.
\end{equation}
We write $[A^0_V(p),X]:= A^0_V(p) (X)$, then it follows by direct computation that the Jacobi identity
\begin{equation}\label{jacobia0}
[A^0_V(p),[S,T]]=[A^0_V(p) (S),T]- [A^0_V(p)(T),S],
\end{equation}
holds for every $S\in\Lie{g}_s$ and $T\in \Lie{g}_t$. This implies that $A^0_V(p)\in {\rm Der}_0(\Lie{g})$, the strata preserving derivations of $\Lie{g}$.
Observe that $A^0_V (p)$ and $-{\rm ad}(V) (p)$ coincide when restricted to the horizontal space, so that conditions~\eqref{conformal} can be read for $A^0_V(p)$. This implies that $A^0_V(p)$ must lie in
$$
\Lie{g}_0=\{D\in {\rm Der}_0(\Lie{g})| D|_{\Lie{g}_{-1}}\in \Lie{co}(2)\}.
$$
The fact that $A^0_V(p)\in\Lie{g}_0$ imposes conditions on higher order derivatives of the coefficients of $V$.

In order to increase the order of derivatives of the coefficients of $V$, we proceed by defining the linear map
$$
A^1_V: \Lie{g}\rightarrow \Lie{g}+\Lie{g}_0
$$
via $A^1_V(p)(X)=\tilde{X}(A^{j+1}_V)(p)$, for every $X\in\Lie{g}_j$.
Namely
\begin{align}\label{A1}
A^1_V(p) (X_1)&=\bmatrix \tilde{X}_1^2f_1(p)&\tilde{X}_1\tilde{X}_2f_1(p)&0&0\\
\tilde{X}_1^2f_2(p)&\tilde{X}_1\tilde{X}_2f_2(p)&0&0 \\
0&0 &\tilde{X}_1\tilde{Y}g(p)&0 \\
0&0&0&\tilde{X}_1\tilde{Z}h(p)
\endbmatrix \nonumber\\
&\nonumber\\
A^1_V(p) (X_2)&=\bmatrix \tilde{X}_2\tilde{X}_1f_1(p)&\tilde{X}_2^2f_1(p)&0&0\\
\tilde{X}_2\tilde{X}_1f_2(p)&\tilde{X}_2^2f_2(p)&0&0\\
0&0 &\tilde{X}_2\tilde{Y}g(p)&0\\
0&0&0&\tilde{X}_2\tilde{Z}h(p)
\endbmatrix \nonumber\\
&\\
A^1_V (p) (Y) &= (\tilde{Y}f_1(p),\tilde{Y}f_2(p),0,0)\nonumber\\
&\nonumber\\
A^1_V (p) (Z) &= (0,0,\tilde{Z}g(p),0).\nonumber
\end{align}
Writing $[A^1_V(p),X]:=A^1_V(p)(X)$ we obtain the Jacobi identity:
$$
[A^1_V(p),[S,T]]=[A^1_V(p) (S),T]- [A^1_V(p)(T),S],
$$
for every $S\in\Lie{g}_s$ and $T\in \Lie{g}_t$, which implies that $A^1_V (p)$ lies in the space
\begin{align*}
\Lie{g}_1 :=\{&u:\Lie{g}\rightarrow \Lie{g}+\Lie{g}_0\;|\;u(\Lie{g}_j)\subset \Lie{g}_{j+1},\\
 &u[S,T]=[u(S),T]-[u(T),S],\forall S\in\Lie{g}_s, \forall T\in\Lie{g}_t\}.
\end{align*}
In general, we continue this procedure inductively and introduce the linear maps
$$A^i_V(p) :\Lie{g}\rightarrow \Lie{g}+\Lie{g}_0 +\Lie{g}_1 +\dots +\Lie{g}_{i-1},$$
 that can vary in a space
\begin{align*}
 \Lie{g}_i :=\{&u:\Lie{g}\rightarrow \Lie{g}+\Lie{g}_0 +\Lie{g}_1 +\dots +\Lie{g}_{i-1}|u(\Lie{g}_j)\subset \Lie{g}_{j+i},\\
& u[S,T]=[u(S),T]-[u(T),S],\forall S\in\Lie{g}_s, \forall T\in\Lie{g}_t\}.
\end{align*}
We say that $\Lie{g}_i$ is the ith {\it prolongation space} of $\Lie{g}$ through $\Lie{g}_0$, and the fact that $A^i_V(p)\in\Lie{g}_i$ provides information on the $i+1$ order derivatives of $f_1$ and $f_2$, the $i+2$ order derivatives of $g$ and $i+3$ order derivatives of $h$.
If the process is finite, then one ends up with a graded Lie algebra $\Lie{g}+\Lie{g}_0+\sum_{k\geq 1} \Lie{g}_k$, where $[u,X]:=u(X)$, for every $X\in\Lie{g}$ and $u\in\sum_{j\geq 0}\Lie{g}_k$.
We shall see that this algebra is isomorphic to the space of vector fields.

\subsection{The Lie algebra of conformal vector fields and the group of conformal maps.}
We proceed by computing the prolongation spaces. First, we compute $\Lie{g}_0$. Write $D= d_{ij} \in \Lie{gl}(4,\R)$. Since $D$ is a strata preserving derivation it follows that
\begin{align*}
 DY&=D[X_1,X_2]=[\sum_{i=1}^2 d_{i1}X_i,X_2] +[X_1,\sum_{j=1}^2 d_{j2}X_j]=(d_{11}+d_{22})Y\\
DZ&= D[X_1,Y]=[\sum_{i=1}^2 d_{i1}X_i,Y]+[X_1,(d_{11}+d_{22})Y]=(2d_{11}+d_{22})Z.
\end{align*}
Since $[X_2,Y]=0$, we have
$$
0= D[X_2,Y]= [\sum_{j=1}^2 d_{j2}X_j,Y]=d_{12}Z,
$$
whence $d_{12}=0$. The condition $D|_{\Lie{g}_{-1}}\in \Lie{co}(2)$ implies $d_{11}=d_{22}$ and $d_{21}=-d_{12}=0$.
In conclusion, $\Lie{g}_0 = \R D$, where $D= {\rm diag}\{1,1,2,3\}$.

The calculation of $\Lie{g}_1$ goes as follows. If $u\in \Lie{g}_1$, then we set $u(X_1)= aD$ and $u(X_2)=bD$ with $a,b\in\R$. By the Jacobi identity we obtain
$u(Y)=aX_2-bX_1$ and $u(Z)=3aY$. Therefore
$$
0=u[X_1,Z]=[aD,Z]-[3aY,X_1]=6aZ,
$$
whence $a=0$, and
$$
0=u[X_2,Z]=[bD,Z] - 0= 3bZ,
$$
whence $b=0$. Thus we conclude that $u=0$ and so $\Lie{g}_1=\{0\}$.

The contact equations \eqref{contact}, the formula \eqref{A0} and the fact that $A_V^0 (p)\in \Lie{g}_0$, lead to the differential  equations
\begin{align}\label{g0}
 \tilde{X}_1f_1 &= \tilde{X}_2 f_2, \qquad \tilde{X}_2 f_1 = -\tilde{X}_1 f_2=0 \nonumber\\
\tilde{Y}g &=2\tilde{X}_1f_1 \\ \tilde{Z}h&=3\tilde{X}_1f_1 \nonumber
\end{align}
Since $A_V^1(p)\in \Lie{g}_1=\{0\}$, equations \eqref{A1} and \eqref{g0} yield
\begin{equation}\label{g1}
 \tilde{X}_1^2 f_1=0 \qquad \tilde{X}_2^2 f_2=0.
\end{equation}
The differential equations above and the contact equations lead to a system of differential equations for all the coefficients of $V$. Moreover, \eqref{contact} implies that the coefficients of $V$ are determined by $h$. So we can restrict our attention to the equations involving $h$:
$$
\tilde{X}_1^3h=0\quad\tilde{X}_2h=0\quad\tilde{Y}^2h=0\quad\tilde{Z}^2h=0.
$$
It is easy to verify that $h$ then varies in a space of polynomials of dimension $5$, so that the space of conformal vector fields, say $\mathcal{C}(G)$, has dimension $5$. In fact, the Tanaka prolongation of $\Lie{g}$ through $\Lie{g}_0$, namely $\Lie{s}=\Lie{g}+\Lie{g}_0$, is isomorphic as Lie algebra to $\mathcal{C}(G)$.
The isomorphism  $\tau: \Lie{s}\rightarrow \mathcal{C}(G)$ is
\begin{equation}\label{tau}
 \tau(X)f (p) =\frac{d}{dt}f({\rm exp}tX \cdot p)|_{t=0},
\end{equation}
where  $\cdot$  denotes the action of ${\rm exp}\Lie{s}$ on $\Lie{g}$. More precisely, the action $\cdot$ is the product on $G$ if $X\in\Lie{g}$, and it indicates the action of the automorphism
${\rm exp}tX$ if $X\in\Lie{g}_0$. Abusing the notation, a basis of $\Lie{s}$ is $\{X_1,X_2,Y,Z,D\}$, and a direct calculation shows
\begin{align*}
\tau(D)&=(3z-2x_1y+\frac{x_1^2x_2}{2})\tilde{Z}+(2y-x_1x_2)\tilde{Y}+x_1\tilde{X}_1+x_2\tilde{X}_2\\
 \tau(X_1)&=(y-x_1x_2)\tilde{Z}+x_2\tilde{Y}+\tilde{X}_1\\
\tau(X_2)&= \frac{x_1^2}{2}\tilde{Z}-x_1\tilde{Y}+\tilde{X}_2\\
\tau(Y)&= -x_1\tilde{Z} +\tilde{Y}\\
\tau(Z)&= \tilde{Z}.
\end{align*}

We showed that the space of vector fields whose local flow is given by conformal mappings is finite dimensional and it coincides with $\Lie{s}$. Now we
prove that if $\phi: \mathcal{U}\rightarrow \mathcal{V}$ is conformal, then $\phi$ is the restriction to $\mathcal{U}$ of the action of some element in ${\rm Aut}\Lie{s}$, the automorphism group of $\Lie{s}$.
The conclusion will be that the space of conformal maps is contained in ${\rm Aut}(\Lie{s})$ and contains ${\rm exp}\Lie{s}$.
By composing with left translations, we may assume that  $e\in\mathcal{U}\cap\mathcal{V}$, and it is enough to show that any conformal map which preserves the identity is such a restriction.

 We show that $\phi$ induces an automorphism of $\Lie{s}$. If $V\in\tau(\Lie{s})$, and we denote by $\psi_t$ the corresponding flow, then $\phi_*V$ is the infinitesimal generator of the 1-parameter group $\phi\psi_t\phi^{-1}$, which is made by conformal maps. Therefore $\phi_*V$ is conformal, whence $\phi_*V\in\tau(\Lie{s})$, and $\tau^{-1}\phi_*\tau\in{\rm Aut}\Lie{s}$.
Since the domain of $\phi$ is connected and contains $e$, the fact that $\phi_*$ is determined by an automorphism of $\Lie{s}$ implies that $\phi$ is also determined by an automorphism of $\Lie{s}$.

It is worth noticing that not all automorphisms of $\Lie{s}$ define conformal maps. For example, let $\alpha$ be the automorphism of $\Lie{s}$ defined by $\alpha(X_1)=X_1$
and $\alpha(X_2)=2X_2$. This automorphism cannot arise as $\tau^{-1}\phi_*\tau$ for some conformal map $\phi$.


\begin{thebibliography}{99}


\bibitem{CC}
{\sc L. Capogna and M. Cowling,}
\newblock \emph{Conformality and Q-harmonicity in Carnot groups}
\newblock{Duke Math. J.}, 135, no. 3, 2006.

\bibitem{cdkr0}
{\sc M. Cowling, F. De Mari, A. Kor\'anyi and H.M. Reimann,}
\newblock \emph{Contact and conformal maps on Iwasawa $N$ groups.}
\newblock Rend. Mat.
Acc. Lincei s.9, vol.\ 13, 2002, 219-232.

\bibitem{CDKR}
{\sc M. Cowling, F. De Mari, A. Kor\'anyi and H.M. Reimann,}
\newblock \emph{ Contact
and conformal mappings in parabolic geometry. I}
\newblock {Geom. Dedicata}, 111:65--86, 2005.	

\bibitem{alefil}
 {\sc F. De Mari and A. Ottazzi,}
 \emph{Rigidity of Carnot groups relative to multicontact structures}, to appear in  Proc. Amer. Math. Soc..

\bibitem{Geh}
{\sc F.W. Gehring,}
\emph{Rings and quasiconformal mappings in space}
\newblock {Trans. Amer. Math. Soc.}, 103:353--393, 1962.

\bibitem{N}
{\sc R. Nevanlinna}
\newblock ``On differentiable mappings'' in \emph{Analytic Functions}, Princeton Math. Ser. 24:3--9, Princeton Univ. Press, Princeton, 1960.


\bibitem{otz1}
{\sc A. Ottazzi,}
\newblock \emph{A sufficient condition for nonrigidity of {C}arnot groups.}
\newblock {Math. Z.}, 259: 617--629, 2008.



\bibitem{otz2}
{\sc A. Ottazzi,}
\newblock \emph{Multicontact vector fields on {H}essenberg manifolds.}
\newblock {J. Lie Theory}, 15(2):357--377, 2005.

\bibitem{OWnote}
{\sc A. Ottazzi and B. Warhurst,}
\newblock \emph{Rigidity of Iwasawa nilpotent Lie groups via Tanaka's theory.}
\newblock Submitted.

\bibitem{OW1}
{\sc A. Ottazzi and B. Warhurst,}
\newblock \emph{Algebraic prolongation and rigidity of Carnot groups}, Monatsh. Math., 2009, DOI 10.1007/s00605-009-0170-7.

\bibitem{Pu90}
{\sc P. Pansu,}
{\em M\'etriques de Carnot-Carath\'eodory et quasiisom\'etries
des espaces sym\'etriques de rang un}.
Ann.\ of Math.\ (2), vol.\ 129, no.\ 1, 1989, 1-60.

		
\bibitem{R}
{\sc H.M. Reimann,}
\newblock \emph{Rigidity of H-type groups}.
\newblock {Math. Z.}, 237 (4): 697-725, 2001.

\bibitem{Res}
{\sc Ju.G. Re\v setnjak [Yu. G. Reshetnyak],}
\emph{Liouville's conformal mapping theorem under minimal regularity hypotheses} (in Russian)
\newblock Sibirsk. Mat. \v Z., 8:835--840, 1967; English translation in Sib. Math. J. 8:631--634, 1967.

\bibitem{ss}
{\sc I. M. Singer and S. Sternberg,}
\emph{The infinite groups of {L}ie and {C}artan. {I}. {T}he
              transitive groups},
\newblock {J. Analyse Math.}, {15:1--114}, {1965}.

\bibitem{tanak1}
{\sc N.~Tanaka,}
\newblock \emph{On differential systems, graded {L}ie algebras and pseudogroups.}
\newblock {J. Math. Kyoto Univ.}, 10: 1--82, 1970.

\bibitem{War1}
{\sc B.~Warhurst,}
\newblock \emph{Jet spaces as nonrigid Carnot groups.}
\newblock J. Lie Theory, 15 (1): 341--356, 2005.

\bibitem{war2}
{\sc B.~Warhurst,}
\newblock \emph{Tanaka prolongation of free Lie algebras.}
\newblock Geom. Dedicata, 130: 59--69, 2007.



\bibitem{yam}
{\sc K.~Yamaguchi,}
\newblock \emph{Differential systems associated with simple graded {L}ie algebras.}
\newblock In Progress in differential geometry, volume~22 of {\em Adv.
  Stud. Pure Math.}, pages 413--494. Math. Soc. Japan, Tokyo, 1993.





\end{thebibliography}
\end{document}